\documentclass[12pt]{article}
\usepackage{epsfig,amsmath,amsthm,amssymb,graphics,amsfonts,natbib}
\usepackage[letterpaper,margin=1in]{geometry}

\theoremstyle{plain}
\newtheorem{thm}{Theorem}[section]
\newtheorem{prop}[thm]{Proposition}

\bibpunct{[}{]}{;}{a}{,}{,}

\begin{document}
 
\date{}
\title{From First Lyapunov Coefficients to Maximal Canards}
\author{Christian Kuehn\thanks{Center for Applied Mathematics, Cornell University}}

\maketitle

\begin{abstract}
Hopf bifurcations in fast-slow systems of ordinary differential equations can be associated with surprising rapid growth of periodic orbits. This process is referred to as canard explosion. The key step in locating a canard explosion is to calculate the location of a special trajectory, called a maximal canard, in parameter space. A first-order asymptotic expansion of this location was found by Krupa and Szmolyan \citep*{KruSzm2,KruSzm1,KruSzm3} in the framework of a ``canard point''-normal-form for systems with one fast and one slow variable. We show how to compute the coefficient in this expansion using the first Lyapunov coefficient at the Hopf bifurcation thereby avoiding use of this normal form. Our results connect the theory of canard explosions with existing numerical software, enabling easier calculations of where canard explosions occur.
\end{abstract}

\section{Introduction}
Our framework in this paper is the theory of fast-slow ordinary differential equations (ODEs):
\begin{eqnarray}
\label{eq:basic1}
\epsilon \dot{x}&=&\epsilon\frac{dx}{d\tau}=f(x,y,\lambda,\epsilon)\\
\dot{y}&=&\frac{dy}{d\tau}=g(x,y,\lambda,\epsilon)\nonumber
\end{eqnarray} 
where $(x,y)\in\mathbb{R}^m\times \mathbb{R}^n$, $\lambda\in\mathbb{R}$ is viewed as a parameter and $\epsilon$ is sufficiently small, i.e. $0<\epsilon\ll 1$. The functions $f:\mathbb{R}^m\times\mathbb{R}^n\times\mathbb{R}\times \mathbb{R}\rightarrow \mathbb{R}^m$ and $g:\mathbb{R}^m\times\mathbb{R}^n\times\mathbb{R} \times \mathbb{R}\rightarrow \mathbb{R}^n$ are assumed to be at least $C^3$ in this paper. The variables $x$ are fast and the variables $y$ are slow. An introduction to the theory of fast-slow systems from the geometric viewpoint can be found in \citep*{ArnoldEncy,Jones,GuckenheimerNDC}, asymptotic methods are developed in \citep*{MisRoz,Grasman} and ideas from nonstandard analysis are considered in \citep*{DienerDiener}. We will only use geometric and asymptotic methods here.\\

In the singular limit $\epsilon\rightarrow 0$ the system \eqref{eq:basic1} becomes a differential-algebraic equation. The algebraic constraint defines the critical manifold:
\begin{equation*}
C_0=\{(x,y)\in\mathbb{R}^m\times \mathbb{R}^n:f(x,y,\lambda,0)=0\}
\end{equation*}
For a point $p\in C_0$ we say that $C_0$ is normally hyperbolic at $p$ if all the  eigenvalues of the $m\times m$ matrix $D_xf(p)$ have non-zero real parts. A normally hyperbolic subset of $C_0$ is an actual manifold and we can locally parametrize it by a map $\psi(y)=x$. This yields the slow subsystem (or reduced flow) $\dot{y}=g(\psi(y),y,\lambda,0)$ defined on $C_0$.\\

Changing in \eqref{eq:basic1} from the slow time scale $\tau$ to the fast time scale $t=\tau/\epsilon$ yields:
\begin{eqnarray}
\label{eq:basic2}
x'&=&\frac{dx}{dt}=f(x,y,\lambda,\epsilon)\\
y'&=&\frac{dy}{dt}=\epsilon g(x,y,\lambda,\epsilon)\nonumber
\end{eqnarray} 
Taking the singular limit $\epsilon\rightarrow 0$ in \eqref{eq:basic2} gives the fast subsystem (or layer equations) $x'=f(x,y,\lambda,0)$ with the slow variables $y$ acting as parameters. A point $p\in C_0$ is an equilibrium of the fast subsystem. We call a subset $S\subset C_0$ an attracting critical manifold if all points $p$ on it are stable equilibria of the fast subsystem i.e. all eigenvalues of $D_xf(p)$ have negative real parts. The subset $S\subset C_0$ is called a repelling critical manifold if for all $p\in S$ at least one eigenvalue of $D_xf(p)$ has positive real part.\\

Fenichel's Theorem \citep*{Fenichel4} states that normally hyperbolic critical manifolds perturb to invariant slow manifolds $C_\epsilon$. A slow manifold $C_\epsilon$ is $O(\epsilon)$ distance away from $C_0$. The flow on the (locally) invariant manifold $C_\epsilon$ converges to the slow subsystem on the critical manifold as $\epsilon\rightarrow 0$. Slow manifolds are usually not unique for a fixed value of $\epsilon=\epsilon_0$ but lie at a distance $O(e^{-k/\epsilon_0})$ away from each other for some $k>0$; nevertheless we shall refer to ``the slow manifold'' associated to subset of the a critical manifold with the possibility of an exponentially small error being understood.\\

Suppose the critical manifold can be divided into two subsets $S_a$ and $S_r$ where $S_a$ is attracting and $S_r$ is repelling so that $C_0=S_a\cup L\cup S_r$. Here $L$ denotes the part of $C_0$ that is not normally hyperbolic. We assume that for $p\in L$ the matrix $D_xf(p)$ has a single zero eigenvalue with right and left eigenvectors $v$ and $w$ and that $w\cdot D_{xx}(p)(v,v)$ and $w\cdot D_yf(p)$ are non-zero. In this case points in $L$ are called fold points. We can use the flow of \eqref{eq:basic2} to extend the associated slow manifolds $S_{a,\epsilon}$ and $S_{r,\epsilon}$ but the extensions might not be normally hyperbolic. The key definition used in this paper is that a trajectory $\gamma$ in the intersection of $S_{a,\epsilon}$ and $S_{r,\epsilon}$ is called a maximal canard; note that this definition requires the extensions of the slow manifolds under the flow. Observe that $\gamma \subset S_{r,\epsilon}$ despite the fact that $S_{r,\epsilon}$ is repelling in the fast directions.\\

We are interested in the case when a fast-slow system undergoes a Hopf bifurcation and a maximal canard is formed close to this bifurcation. The periodic orbits resulting from the Hopf bifurcation grow rapidly in a $\lambda$-interval of width $O(e^{-K/\epsilon})$ for some $k>0$. The rapid orbit growth is usually referred to as canard explosion and the bifurcation scenario is called singular Hopf bifurcation.\\

The paper is organized as follows. In Section \ref{sec:ce} we describe results on singular Hopf bifurcation and canard explosion obtained by Krupa and Szmolyan \citep*{KruSzm2}. In Section \ref{sec:Hopf} we clarify the different definitions of the first Lyapunov coefficient of a Hopf bifurcation. In Section \ref{sec:relation} we present the main results on the relation between the location of the maximal canard and the first Lyapunov coefficient. We describe which terms will contribute to the first order approximation using a rescaled Hopf bifurcation normal form. Then we show explicitly how to compute a first order approximation to the location of the maximal canard avoiding additional center manifold reduction and normal form transformations. In Section \ref{sec:examples} we locate the maximal canards in two examples: a two-dimensional version of van der Pol's equation and a three-dimensional version of the FitzHugh-Nagumo equation.\\

Note that we do not give a detailed description of dynamics associated to a singular Hopf bifurcation and refer the reader to the previous extensive literature e.g. \citep*{BaerErneuxI,BaerErneuxII,Braaksma,KruSzm2,GuckenheimerSH}. 
 
\section{Canard Explosion}
\label{sec:ce}

We describe the main results about canard explosion in fast-slow systems with one fast and one slow variable from \citep*{KruSzm2}. Consider a planar fast-slow system of the form
\begin{eqnarray}
\label{eq:fsgen}
x'&=& f(x,y,\lambda,\epsilon)\nonumber\\
y'&=& \epsilon g(x,y,\lambda,\epsilon)
\end{eqnarray}
where $f,g\in C^k(\mathbb{R}^4,\mathbb{R})$ for $k\geq 3$, $\lambda\in\mathbb{R}$ is a parameter and $0<\epsilon\ll 1$. Denote the critical manifold of \eqref{eq:fsgen} by $C_0$. We assume that $C_0$ is locally parabolic with a minimum at the origin $(x,y)=(0,0)$ independent of $\lambda$ so that $(0,0)$ is a fold point; more precisely 
\begin{equation}
\label{eq:fold}
f(0,0,\lambda,0)=0, \quad f_x(0,0,\lambda,0)=0,\quad f_{xx}(0,0,\lambda,0)\neq 0,\quad f_{y}(0,0,\lambda,0)\neq 0
\end{equation}
In addition, we assume that $g(0,0,\lambda\neq 0,0)\neq 0$; under these conditions the fold point at the origin is generic for $\lambda\neq 0$. We assume without loss of generality that $f_{xx}(0,0,\lambda,0)>0$ so that $C_0$ is locally a parabola with a minimum at the origin.  Using \eqref{eq:fold} and the implicit function theorem we have that $C_0$ is the graph of a function $y=\phi(x)$ for $\phi:U\rightarrow \mathbb{R}$ where $U$ is a sufficiently small neighbourhood of $x=0$. Assume that $C_0$ splits into an attracting and a repelling curve $C=C_l\cup \{(0,0)\}\cup C_r$ where
\begin{equation*}
C_l=\{x<0,f_x<0\}\cap C_0,\quad C_r=\{x>0,f_x>0\}\cap C_0
\end{equation*}

\begin{figure}[htbp]
\centering
 \includegraphics[width=0.96\textwidth]{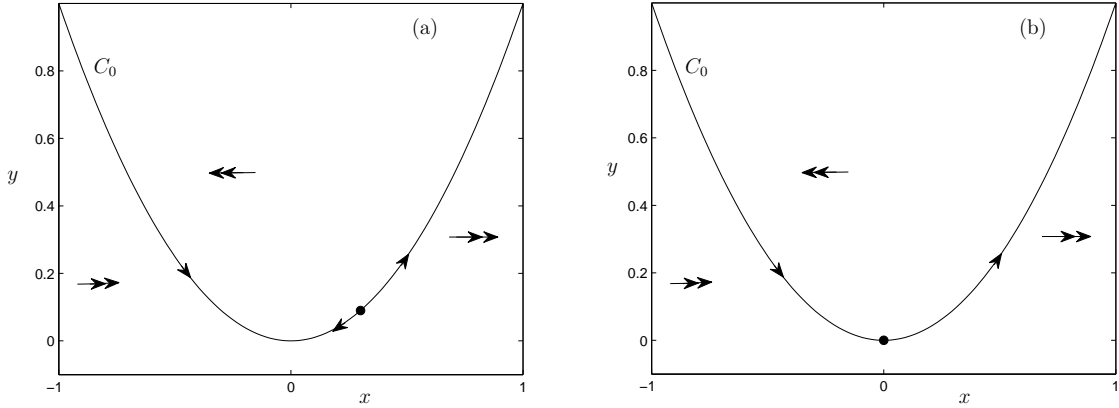} 
\caption{\label{fig:fold}(a) A generic fold for $\lambda\neq 0$. (b) A nondegenerate canard point for $\lambda=0$. The slow flow is indicated by single and the fast flow by double arrows.}
\end{figure}

The situation is shown in Figure \ref{fig:fold}(a). Differentiating $y=\phi(x)$ with respect to $\tau=t\epsilon$ we get that the slow flow on $C_0$ is defined by 
\begin{equation*}
\frac{dx}{d\tau}=\dot{x}=\frac{g(x,\phi(x),\lambda,0)}{\phi'(x)}
\end{equation*}
Note that the slow flow is singular for $\lambda \neq 0$ at $(0,0)$ since $\phi'(0)=0$ and $g(0,0,\lambda,0)\neq0$. Assume that at $\lambda=0$ we have a non-degenerate canard point (see Figure \ref{fig:fold}(b)) so that in addition to the fold conditions we have
\begin{equation*}
g(0,0,0,0)=0,\quad g_x(0,0,0,0)\neq 0,\quad g_\lambda(0,0,0,0)\neq 0
\end{equation*}
Therefore the slow flow is well-defined at $(0,0)$ for $\lambda=0$ and we assume without loss of generality that $\dot{x}>0$ in this case. Near a non-degenerate canard point \eqref{eq:fsgen} can be transformed into a normal form \citep*{KruSzm3}:
\begin{eqnarray}
\label{eq:nform}
x'&=& -yh_1(x,y,\lambda,\epsilon)+x^2h_(x,y,\lambda,\epsilon)+\epsilon h_3(x,y,\lambda,\epsilon)\nonumber\\
y'&=& \epsilon(xh_4(x,y,\lambda,\epsilon)-\lambda h_5(x,y,\lambda,\epsilon)+yh_6(x,y,\lambda,\epsilon))
\end{eqnarray}
where the functions $h_i$ are given by:
\begin{eqnarray*}
h_3(x,y,\lambda,\epsilon)&=&O(x,y,\lambda,\epsilon)\\
h_j(x,y,\lambda,\epsilon)&=&1+O(x,y,\lambda,\epsilon),\qquad j=1,2,4,5
\end{eqnarray*}
We define several computable constants, abbreviating $(0,0,0,0)=0$ in the definitions:
\begin{equation*}
a_1=(h_3)_x(0),\quad a_2=(h_1)_x(0),\quad a_3=(h_2)_x(0),\quad a_4=(h_4)_x(0),\quad a_5=(h_6)_x(0) 
\end{equation*}
Note that all $a_i$ for $i=1,2,3,4,5$ only depend on partial derivatives with respect to $x$. Next we define another constant:
\begin{equation*}
A=-a_2+3a_3-2a_4-2a_5
\end{equation*}

\begin{thm}(\citep*{KruSzm2}) 
\label{thm:main} For $0<\epsilon<\epsilon_0$, $|\lambda|<\lambda_0$ and $\epsilon_0>0$, $\lambda_0>0$ sufficiently small and under the previous assumptions in this section there exists a unique equilibrium point $p$ for $\eqref{eq:nform}$ in a neighbourhood of $(x,y)=(0,0)$. The equilibrium $p$ undergoes a Hopf bifurcation at $\lambda_H$ with
\begin{equation}
\label{eq:hopf}
\lambda_H=-\frac{a_1+a_5}{2}\epsilon+O(\epsilon^{3/2})
\end{equation}
The slow manifolds $C_{\epsilon,l}$ and $C_{\epsilon,r}$ intersect/coincide in a maximal canard at $\lambda_c$ for
\begin{equation}
\label{eq:canard}
\lambda_c=-\left(\frac{a_1+a_5}{2}+\frac{A}{8}\right)\epsilon+O(\epsilon^{3/2})
\end{equation} 
The equilibrium $p$ is stable for $\lambda<\lambda_H$ and unstable for $\lambda>\lambda_H$. The Hopf bifurcation is non-degenerate for $A\neq 0$, supercritical for $A<0$ and subcritical for $A>0$.
\end{thm}

\textit{Remark:} The asymptotic expansions for $\lambda_H$ and $\lambda_c$ are asymptotic series with asymptotic sequence $\{\epsilon^{k/2}\}_{k=0}^{\infty}$ and Theorem \ref{thm:main} implies that the first two coefficients of the expansion are zero and the third coefficient can be computed explicitly.\\

Note that currently no standard bifurcation software such as AUTO \citep*{Doedel_AUTO2007} or MatCont \citep*{MatCont} computes the constants $a_i$ and $A$ automatically. Nevertheless bifurcation software can detect Hopf bifurcations so that given a fixed $\epsilon$ we can approximate $\lambda_H$ numerically. Hence the numerical problem that remains is to compute A since
\begin{equation*}
\lambda_H-\lambda_c=\frac{A}{8}\epsilon+O(\epsilon^{3/2})
\end{equation*}  
To simplify the notation we define $K=A/8$. If we know $K$ we can easily approximate the location of the maximal canard by $\lambda_c=\lambda_H-K\epsilon+O(\epsilon^{3/2})$. The maximal canard organizes the canard explosion \citep*{KruSzm2} and indicates where the rapid amplitude growth of the small orbits generated in the Hopf bifurcation occurs. Our goal is to avoid any additional normal form transformations and center manifold reductions to compute $K$. The key point to achieve this is to observe that $K$ is just a rescaled version of ``the'' first Lyapunov coefficient of the Hopf bifurcation at $\lambda_H$. 

\section{The First Lyapunov Coefficient}
\label{sec:Hopf}

We review and clarify the interpretation, computation and conventions associated with the first Lyapunov coefficient of a Hopf bifurcation. Consider a general N-dimensional ODE at a non-degenerate Hopf bifurcation point. We assume that the equilibrium has been translated to the origin so that
\begin{equation}
\label{eq:hopfnf}
z'=Mz+F(z),\qquad \text{for $z\in\mathbb{R}^N$} 
\end{equation} 
with $F(z)=O(\|z\|^2)$ and $M=(m_{ij})$. Taylor expanding $F$ yields
\begin{equation*}
z'=Mz+\frac12 B(z,z)+\frac16C(z,z,z)
\end{equation*} 
where the multilinear functions $B$ and $C$ are given by:
\begin{eqnarray*}
B_i(u,v)&=&\sum_{j,k=1}^N \left.\frac{\partial^2 F_i(\xi)}{\partial \xi_j\partial\xi_k}\right|_{\xi=0}u_jv_k\\
C_i(u,v,w)&=&\sum_{j,k,l=1}^N \left.\frac{\partial^3 F_i(\xi)}{\partial \xi_j\partial\xi_k \partial \xi_l}\right|_{\xi=0}u_jv_kw_l
\end{eqnarray*}
The matrix $M$ has eigenvalues $\lambda_{1,2}=\pm i\omega_0$ for $\omega_0>0$. Let $q\in \mathbb{C}^N$ be the eigenvector of $\lambda_1$ and $p\in\mathbb{C}^N$ the corresponding eigenvector of the transpose $M^T$ i.e.
\begin{eqnarray*}
Mq=i\omega_0 q, \quad M\bar{q}=-i\omega_0\bar{q},\quad M^Tp=-i\omega_0 p, \quad M\bar{p}=i\omega_0\bar{p} 
\end{eqnarray*}
where the the overbar denotes componentwise complex conjugation. We can always normalize $p$ so that the standard complex inner product with $q$ satisfies $\bar{p}^T q=\sum_{j=1}^N\bar{p}_jq_j=1$. The first Lyapunov coefficient of the Hopf bifurcation can then be defined by (\citep*{Kuznetsov}, p.180):
\begin{equation}
\label{eq:lyaK1}
l^{Ku}_1=\frac{1}{2\omega_0}\left(\bar{p}^T C(q,q,\bar{q})-2\bar{p}^T B(q,L^{-1}B(q,\bar{q}))+\bar{p}^T B(\bar{q},(2i\omega_0I_N-M)^{-1}B(q,q))\right)
\end{equation}
In the case of a two-dimensional vector field $F=(F^1,F^2)$ the formula \eqref{eq:lyaK1} can be expressed in the simpler form (\citep*{Kuznetsov}, p.98):
\begin{equation}
\label{eq:lyaK2}
l_1^{Ku}=\frac{1}{2\omega_0^2}\text{Re}(ig_{20}g_{11}+\omega_0g_{21})
\end{equation} 
where 
\begin{equation*}
g_{20}=\bar{p}^T B(q,q),\quad g_{11}=\bar{p}^T B(q,\bar{q}), \quad g_{21}=\bar{p}^T C(q,q,\bar{q})
\end{equation*}
It is important to note that $l^{Ku}_1$ is not uniquely defined until we choose a normalization of the eigenvector $q$. We adopt the convention using unit norm $\bar{q}^Tq=1$. A slight modification of the formula \eqref{eq:lyaK1} is used to evaluate the Lyapunov coefficient $l_1^{MC}$ numerically in the bifurcation software MatCont \citep*{MatCont}. Using the current MatCont convention\footnote{MatCont version 2.5.1 - December 2008} we note that 
\begin{equation*}
\omega_0 l_1^{Ku}=l_1^{MC}
\end{equation*}
Other expressions for the first Lyapunov coefficient can be found in the literature. We consider only the planar case using simpler notation $(z_1,z_2)=(x,y)$: 
\begin{equation}
\label{eq:2D}
\left(\begin{array}{c} x' \\ y'\end{array}\right)= \left(\begin{array}{c} F^1(x,y) \\ F^2(x,y) \end{array}\right)=: M\left(\begin{array}{c} x \\ y \end{array}\right) + \left(\begin{array}{c} f(x,y) \\ g(x,y) \end{array}\right)
\end{equation}
Then another convention for $l_1$ is (\citep*{ChowLiWang}, p.211):
\begin{eqnarray}
\label{eq:lyaC}
l_1^{CLW}=&&\frac{m_{12}}{16 \omega_0^4} [\omega_0^2[(f_{xxx}+g_{xxy}) + 2 m_{22}(f_{xxy}+g_{xyy})-m_{21}(f_{xyy}+g_{yyy})]\nonumber\\
&& -m_{12}m_{22}(f_{xx}^2-f_{xx}g_{xy}-f_{xy}g_{xx}-g_{xx}g_{yy}-2g_{xy})\nonumber\\
&&-m_{21}m_{22}(g_{yy}^2-g_{yy}f_{xy}-g_{xy}f_{yy}-f_{xx}f_{yy}-2f_{xy}^2)\\
&&+m_{12}^2(f_{xx}g_{xx}+g_{xx}g_{xy})-m_{21}^2(f_{yy}g_{yy}+f_{xy}f_{yy})\nonumber\\
&&-(\omega_0^2+3m_{22}^2)(f_{xx}f_{xy}-g_{xy}g_{yy})]\nonumber
\end{eqnarray}
where all evaluations in \eqref{eq:lyaC} are at $(x,y)=(0,0)$. Next, assume that we have applied a preliminary linear coordinate change
\begin{equation*} 
\left(\begin{array}{c} x \\ y\end{array}\right)=N\left(\begin{array}{c} u \\ v\end{array}\right) \qquad \text{where $N:\mathbb{R}^2\rightarrow \mathbb{R}^2$} 
\end{equation*}
to the system \eqref{eq:2D} to transform $M$ into Jordan normal form. Then we look at:
\begin{eqnarray}
\label{eq:2Dt}
\left(\begin{array}{c} u' \\ v'\end{array}\right)&=&  \left(\begin{array}{cc} 0 & -\omega_0 \\ \omega_0 & 0 \end{array}\right) \left(\begin{array}{c} u \\ v \end{array}\right) + N^{-1}\left(\begin{array}{c} f(u,v) \\ g(u,v) \end{array}\right)\nonumber\\
&=&\left(\begin{array}{cc} 0 & -\omega_0 \\ \omega_0 & 0 \end{array}\right) \left(\begin{array}{c} u \\ v \end{array}\right) + \left(\begin{array}{c} f^*(u,v) \\ g^*(u,v) \end{array}\right)
\end{eqnarray}
In this case the Lyapunov coefficient formula simplifies (\citep*{GH}, p.152):
\begin{eqnarray}
\label{eq:lyaGH}
l_1^{GH}=&&\frac{1}{16}[f^*_{xxx}+f^*_{xyy}+g^*_{xxy}+g^*_{yyy}]+\frac{1}{16\omega_0}[f^*_{xy}(f^*_{xx}+f^*_{yy})\nonumber\\
&&-g^*_{xy}(g^*_{xx}+g^*_{yy})-f^*_{xx}g^*_{xx}+f^*_{yy}g^*_{yy}]
\end{eqnarray}
Note that the linear transformation $N$ is not unique. We adopt the convention that
\begin{equation*}
N=\left(\begin{array}{cc} 2\text{Re}(q_1) & -2\text{Im}(q_1) \\ 2\text{Re}(q_2) & -2\text{Im}(q_2) \end{array}\right) 
\end{equation*}
where $q=(q_1,q_2)$ is the normalized eigenvector of the linearization $L$ that satisfies $Lq=i\omega_0 q$. Another common definition for \eqref{eq:2Dt} is (\citep*{Perko}, p.353):
\begin{eqnarray}
\label{eq:lyaP}
l_1^{Pe}=&&\frac{3\pi}{4\omega_0^2}([f^*_{xy}f^*_{yy}+f^*_{yy}g^*_{yy}-f^*_{xx}g^*_{xx}-g^*_{xy}g^*_{xx}-g^*_{xy}g^*_{yy}+f^*_{xy}f^*_{xx}]\nonumber\\
&&+\omega_0[g^*_{yyy}+f^*_{xxx}+f^*_{xyy}+g^*_{xxy}])
\end{eqnarray}
The Hopf bifurcation theorem holds for any version of $l_1$ as only the sign is relevant in this case:
\begin{thm}(see e.g. \citep*{GH,Kuznetsov})
\label{thm:hopf} 
A non-degenerate Hopf bifurcation of \eqref{eq:hopfnf} is supercritical if $l_1<0$ and subcritical if $l_1>0$. 
\end{thm} 
 
Since we need not only a qualitative result such as Theorem \ref{thm:hopf}, but a quantitative one relating the Lyapunov coefficient to canard explosion, it is necessary to distinguish between the different conventions we reviewed above.

\section{Relating $l_1$ and $K$}
\label{sec:relation}

Krupa and Szmolyan consider a blow-up \citep*{KruSzm2,KruSzm3,KruSzm1} of \eqref{eq:nform} given by $\Phi:S^3\times I\rightarrow \mathbb{R}^4$ in a particular chart $K_2$:
\begin{equation}
\label{eq:blowup}
x=rx_2,\qquad y=r^2y_2,\qquad \lambda=r\lambda_2,\qquad  \epsilon=r^2\epsilon_2
\end{equation}
where $r\in I\subset \mathbb{R}$ and $(x_2,y_2,\lambda_2,\epsilon_2)\in S^3$. Using \eqref{eq:blowup} the resulting vector field can be desingularized by dividing it by $\sqrt{\epsilon}$. We shall not discuss the details of the blow-up approach and just note that this transformation and the following desingularization are simply a rescaling of the vector field given by:
\begin{equation}
\label{eq:rescale}
x_2=\epsilon^{-1/2}x,\qquad y_2=\epsilon^{-1}y,\qquad \lambda_2=\epsilon^{-1/2}\lambda, \qquad t_2=\epsilon^{1/2}t
\end{equation}
Using the formula from Chow, Li and Wang \citep*{ChowLiWang} in the rescaled version of \eqref{eq:nform} Krupa and Szmolyan get the following result:
\begin{prop}
\label{prop:bu}
In the coordinates \eqref{eq:rescale} the first Lyapunov coefficient $l_1$ has asympototic expansion:
\begin{equation}
\label{eq:lyaKS}
\bar{l}^{CLW}_1=K\sqrt{\epsilon}+O(\epsilon)
\end{equation}
where the overbar indicates the first Lyapunov coefficient in coordinates given by \eqref{eq:rescale}.
\end{prop}

First, we want to explain in more detail which terms in the vector field \eqref{eq:nform} contribute to the leading order coefficient $K$. The main problem is that the Lyapunov coefficient is often calculated after an $\epsilon$-dependent rescaling, such as \eqref{eq:rescale}, has been carried out. This can lead to rather unexpected effects in which terms contribute to the Lyapunov coefficient, as pointed out by Guckenheimer \citep*{GuckenheimerSH} in the context of singular Hopf bifurcation in $\mathbb{R}^3$. \\

To understand how the rescaling \eqref{eq:rescale} affects the Lyapunov coefficient we consider the Hopf normal form case. We start with a planar vector field with linear part in Jordan form \eqref{eq:2Dt}. Assume that the equilibrium is at the origin $(x,y)=0$ and Hopf bifurcation occurs for $\lambda=0$. Applying the rescaling \eqref{eq:rescale} we get:
\begin{equation}
\label{eq:rescaledGH}
\left(\begin{array}{c} dx_2/dt_2 \\ dy_2/dt_2\end{array}\right)= \left(\begin{array}{cc} 0 & -\omega_0/\sqrt{\epsilon} \\ \omega_0/\sqrt{\epsilon} & 0 \end{array}\right) \left(\begin{array}{c} x_2 \\ y_2 \end{array}\right) + \left(\begin{array}{c} \frac{1}{\epsilon} f^*(\sqrt{\epsilon}x_2,\epsilon y_2) \\ \frac{1}{\epsilon^{3/2}}g^*(\sqrt{\epsilon}x_2,\epsilon y_2) \end{array}\right)
\end{equation}
In a fast-slow system with singular Hoof bifurcation we know that $g^*(.,.)=\epsilon(\ldots)$ and that $\omega_0=O(\sqrt{\epsilon})$. Setting $k_\omega=\omega_0/\sqrt{\epsilon}$ the Lyapunov coefficient can be computed to leading order by \eqref{eq:lyaGH}:
\begin{equation}
\label{eq:GHbu}
\bar{l}_1^{GH}=\frac{1}{k_\omega}\left(f^*_{x_2x_2}(0,0)[g^*_{x_2x_2}(0,0)+f^*_{x_2y_2}(0,0)]+k_\omega f_{x_2x_2x_2}(0,0)\right) \sqrt{\epsilon}+O(\epsilon)
\end{equation}
Equation \eqref{eq:GHbu} explains the leading-order behaviour more clearly and shows that due to the rescaling certain derivative terms in the Lyapunov coefficient for a singular Hopf bifurcation are non-leading terms with respect to $\epsilon\rightarrow 0$. The point is that the rescaling modifies the order with respect to $\epsilon$ of the linear and nonlinear terms. Also, applying the chain rule to the nonlinear terms to calculate the necessary derivatives can affect which terms contribute.\\

To make Proposition \eqref{prop:bu} more useful in an applied framework we have computed all the different versions of the Lyapunov coefficient defined in Section \eqref{sec:Hopf} up to leading order for equation \eqref{eq:nform} in original non-rescaled coordinates. The computer algebra system Maple \citep*{Maple} was used in this case:
\begin{eqnarray}
\label{eq:result}
l_1^{Ku}&=&\frac{4K}{\sqrt{\epsilon}}+O(\sqrt{\epsilon})\nonumber\\
l_1^{MC}&=&\frac{4K\omega_0}{\sqrt{\epsilon}}+O(\omega_0\sqrt{\epsilon})\nonumber\\
l_1^{GH}&=&K+O(\epsilon)\\
l_1^{CLW}&=&K+O(\epsilon)\nonumber\\
l_1^{Pe}&=&\frac{3\pi K}{64\omega_0}+O(\epsilon/\omega_0)\nonumber
\end{eqnarray}
Using the results \eqref{eq:result} we now have a direct strategy how to analyze a canard explosion generated in a singular Hopf bifurcation.
\begin{enumerate}
 \item Compute the location of the Hopf bifurcation. This gives $\lambda_H$.
 \item Find the first Lyapunov coefficient at the Hopf bifurcation, e.g. we get $l_1^{MC}\approx 4K\omega_0/\sqrt{\epsilon}$.
 \item Compute the location of the maximal canard, and hence the canard explosion, by $\lambda_c\approx \lambda_H-K\epsilon$. For example, using MatCont we would get 
\begin{equation}
\label{eq:MCpractical}
\lambda_c\approx \lambda_H-\frac{l_1^{MC}}{4\omega_0}\epsilon^{3/2}
\end{equation}
\end{enumerate}
Observe that the previous calculation may require calculating the eigenvalues at the Hopf bifurcation to determine $\omega_0$ but does not require any center manifold calculations nor additional normal form transformations; these have basically been encoded in the calculation of the Lyapunov coefficient.

\section{Examples}
\label{sec:examples}

\begin{figure}[htbp]
\centering
 \includegraphics[width=0.8\textwidth]{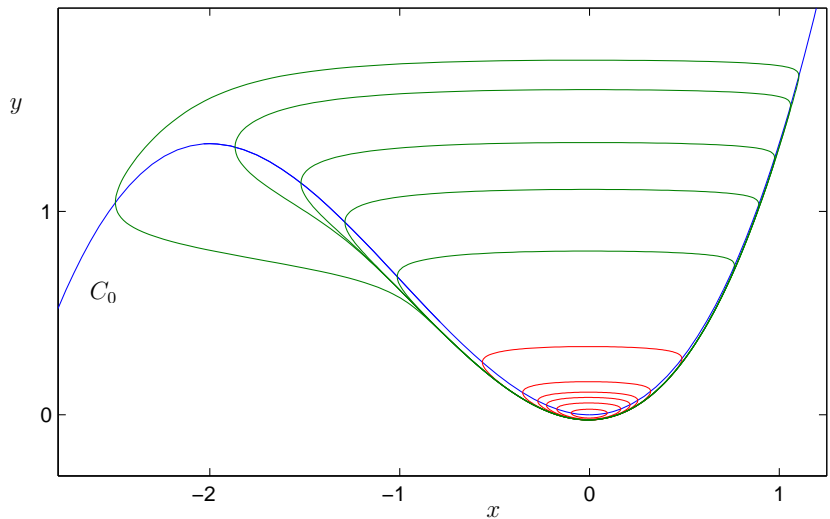} 
\caption{\label{fig:vdP}Continuation of periodic orbits emanating from the Hopf bifurcation at $\lambda=0$. The parameter values for the red orbits are $\lambda=-0.001$, $-0.0025$, $-0.004$, $-0.005$, $-0.006$, $-0.0065$ and for all the green orbits the parameter value is $\lambda\approx -0.006509$ indicating a canard explosion near this parameter value.}
\end{figure}

The first example is a version of van der Pol's equation \citep*{vanderPol,vanderPol_RO,KruSzm2} given by:
\begin{eqnarray*}
x'&=&y-x^2-\frac{x^3}{3}\nonumber\\
y'&=&\epsilon(\lambda-x)
\end{eqnarray*} 
We have to reverse time $t\rightarrow -t$ to satisfy the assumptions of Section \eqref{sec:ce}. This gives:
\begin{eqnarray}
\label{eq:vdP}
x'&=&x^2+\frac{x^3}{3}-y\nonumber\\
y'&=&\epsilon(x-\lambda)
\end{eqnarray} 
The critical manifold is given by $C_0=\{y=x^2+x^3/3\}$ with two fold points at $(0,0)$ and $(-2,4/3)$. The fold points split the critical manifold into three normally hyperbolic parts:
\begin{equation*}
C_l=C_0\cap \{x<-2\},\quad C_m=C_0\cap \{-2<x<0\},\quad C_r=C_0\cap \{0<x\}
\end{equation*}
We only study the fold point at the origin which becomes a canard point for $\lambda=0$. The unique equilibrium point $p=(x_e(\lambda),y_e(\lambda))$ of \eqref{eq:vdP} lies on $C_0$ and satisfies $x_e(\lambda)=\lambda$. It is easy to check that subcritical Hopf bifurcation occurs for $\lambda=\lambda_H=0$. Matching terms in \eqref{eq:vdP} and the normal form \eqref{eq:nform} we find:
\begin{equation*}
h_1=h_4=h_5=1,\qquad h_2=1+\frac13x,\quad h_6=0
\end{equation*}
Therefore $K=1/8$ and we find analytically that the location of the maximal canard representing the intersection of $C_{m,\epsilon}$ and $C_{r,\epsilon}$ is 
\begin{equation}
\label{eq:mcvdP}
\lambda_c=-(1/8)\epsilon+O(\epsilon^{3/2})
\end{equation}
A numerical continuation calculation using a bifurcation software tool - we used MatCont \citep*{MatCont} - gives that the first Lyapunov coefficient for the Hopf bifurcation at $\lambda=\lambda_H=0$ for $\epsilon=0.05$ is
\begin{equation*}
l_1^{MC}\approx 0.4762
\end{equation*}  
An easy calculation\footnote{Using MatCont 2.5.1. we can modify the file \texttt{/matcont2.5.1/MultilinearForms/nf\_H.m} to return the variable \texttt{omega}$=\omega_0$ or to return $l_1^K=l_1^{MC}/\omega_0$.} yields that $\omega_0\approx 0.2236$. Using \eqref{eq:MCpractical} we compare this to the result in equation \eqref{eq:mcvdP} with $\epsilon=0.05$. Dropping higher-order terms we have:
\begin{equation}
\label{eq:test1}
\lambda_c(\text{analytical})=-0.0063,\qquad \lambda_c(\text{numerical using $l_1$})=-0.0060
\end{equation}
The coincidence of the values of the location of the maximal canard is already quite good but this is expected since we have only compared the asymptotic formula to the Lyapunov coefficient formula derived from it which was evaluated numerically using continuation. A simple direct test to compare \eqref{eq:test1} to the location of the maximal canard is to use continuation of periodic orbits from the Hopf bifurcation point. The results are shown in Figure \ref{fig:vdP}.\\

\textit{Remark:} Depending on the bifurcation software used, direct continuation of periodic orbits can fail for small values of $\epsilon$. In this case special methods are needed to continue periodic orbits having canard segments; see e.g. \citep*{GuckLaMar,GuckCK2,DesrochesKrauskopfOsinga2}. Note that locating Hopf bifurcations and calculating Lyapunov coefficients works well even for very small values of $\epsilon$ as we require only local algebraic calculations.\\

We conclude from Figure \ref{fig:vdP} that our estimates in \eqref{eq:test1} are very good indicators to determine where the canard explosion exists since they are already decent for a relatively large $\epsilon=0.05$. In many standard fast-slow systems values of $\epsilon\leq 0.01$ are commonly considered.\\

In higher dimensions the analytical calculations will be very difficult to carry out. As a second example consider a version of the FitzHugh-Nagumo equation \citep*{Sneydetal,GuckCK1,GuckCK3}:
\begin{eqnarray}
\label{eq:fhn}
x_1'&=&x_2\nonumber\\
x_2'&=&\frac15\left(sx_2-x_1(x_1-1)(0.1-x_1)+y-I\right)\\
y'&=& \frac{\epsilon}{s}(x_1-y)\nonumber
\end{eqnarray} 
where $I$ and $s$ are parameters. Equation \eqref{eq:fhn} has two fast and one slow variable and a unique equilibrium point $p(I)=p$. For a detailed fast-slow system analysis describing the bifurcations we refer the reader to \citep*{GuckCK1,GuckCK3}. We only note that the critical manifold is cubic curve given by
\begin{equation*}
C_0=\{(x_1,x_2,y)\in\mathbb{R}^3:x_2=0\text{ and }y=x_1(x_1-1)(0.1-x_1)+I\}
\end{equation*}
It is normally hyperbolic away from two fold points $x_{1,\pm}$ given by the local minimum and maximum of the cubic. The equilibrium $p$ passes through the fold points under parameter variation; $O(\epsilon)$ away from these points the equilibrium undergoes Hopf bifurcation \citep*{GuckCK1,GuckCK3}. We shall just compute a particular case applying our result \eqref{eq:result}. We fix $s=1.37$ and observe that $I$ plays the same role as $\lambda$ in our previous calculations. We calculated the location of the maximal canard for several values of $\epsilon$. The results are shown in Table \ref{tab:dist}.\\

\begin{table}[htbp]
\begin{center}
\begin{tabular}{|c|c|c|}
\hline
$\epsilon$       & $I_c$               & $I_c(Lyapunov)$\\
\hline  
$10^{-2}$        & $\approx 0.0582046$ & $\approx 0.06308$ \\
$5\cdot 10^{-3}$ & $\approx 0.0545535$ & $\approx 0.05629$ \\
$10^{-3}$        & $\approx 0.0517585$ & $\approx 0.05196$ \\
$5\cdot 10^{-4}$ & $\approx 0.0514108$ & $\approx 0.05150$ \\
\hline
\end{tabular}
\caption{Comparison between the actual location of the maximal canard (canard explosion) $I_c$ and the first-order approximation $I_c(Lyapunov)$ computed using the first Lyapunov coefficient at the Hopf bifurcation.\label{tab:dist}}
\end{center}
\end{table}

The second column of Table \ref{tab:dist} shows the actual location $I_c$ of the maximal canard (canard explosion) obtained from continuation of periodic orbits using AUTO \citep*{Doedel_AUTO2007}. The third column shows the approximation obtained by using the first Lyapunov coefficient. The Lyapunov coefficient has been computed using MatCont \citep*{MatCont}. The error $E(\epsilon)$ of this calculation is of order $O(\epsilon^{3/2})$ as expected from \eqref{eq:lyaKS}. Obviously the approximation improves for smaller values of $\epsilon$. \\

In the case of $\epsilon=0.01$ it has been shown in \citep*{GuckCK1,GuckCK3} that there is an intricate bifurcation scenario involving homoclinic orbits in a parameter interval near $I_c$. The first order approximation of the maximal canard is not sufficient to relate it to the homoclinic bifurcation. This shows that the magnitude of $\epsilon$ and other relevant bifurcations in the system have to be taken into account carefully when applying the results we presented here.  

\section{Discussion}

We have investigated the relation between the first Lyapunov coefficient at a singular Hopf bifurcation and the associated maximal canard orbit. The major result is that no additional algorithms are needed to compute a first order approximation to the location of the maximal canard. Standard bifurcation software packages compute the Lyapunov coefficient and our results can be used to approximate the maximal canard location from this numerical calculation.\\

We also pointed out that there is no ``standard definition'' of the first Lyapunov coefficient of a Hopf bifurcation. This is not surprising since classical qualitative bifurcation theory only requires the sign of the Lyapunov coefficient. We hope that the comparison in Section \ref{sec:Hopf} will help the reader to adapt their own numerical algorithms and software packages to support the calculation of maximal canard locations.\\

Open questions which we leave for future work include the extensions to multiple slow variables, higher-order asymptotic expansions and the relation between the Lyapunov coefficient and blow-up transformations.

\end{document}